\documentclass[a4paper,12pt]{article}
\makeatletter
\@addtoreset{footnote}{page}
\makeatother

\usepackage{amsthm}
\usepackage{amsmath,amssymb,latexsym,amsfonts,mathrsfs}

\theoremstyle{plain}
\newtheorem{Thm}{Theorem}[section]
\newtheorem{Lem}[Thm]{Lemma}
\newtheorem{Prop}[Thm]{Proposition}
\newtheorem{Cor}[Thm]{Corollary}
\theoremstyle{definition}

\newtheorem{Rem}[Thm]{Remark}

\newcommand{\Proof}[2][Proof]{\begin{proof}[{#1}] #2 \end{proof}}

\numberwithin{equation}{section}

\makeatletter
\renewcommand\section{\@startsection {section}{1}{\z@}%
                                   {-3.5ex \@plus -1ex \@minus -.2ex}%
                                   {2.3ex \@plus.2ex}%
                                   {\normalfont\large\bf}}
\makeatother

\makeatletter
\renewcommand\subsection{\@startsection {subsection}{1}{\z@}%
                                   {-3.5ex \@plus -1ex \@minus -.2ex}%
                                   {2.3ex \@plus.2ex}%
                                   {\normalfont\normalsize\bf}}
\makeatother

\renewcommand{\d}{{\rm d}} %differential
\newcommand{\eps}{\ensuremath{\varepsilon}}

\newcommand{\dist}{\stackrel{{\rm d}}{=}}
\newcommand{\tend}[2]{\mathrel{\mathop{\longrightarrow}\limits^{#1}_{#2}}}

\renewcommand{\tilde}{\widetilde}
\renewcommand{\bar}{\overline}

\newcommand{\norm}[1]{\left\| #1 \right\|} %norm
\newcommand{\cbra}[1]{\!\left\{ #1 \right\}} %curly brackets or braces
\newcommand{\cA}{\ensuremath{\mathcal{A}}}
\newcommand{\cB}{\ensuremath{\mathcal{B}}}
\newcommand{\cI}{\ensuremath{\mathcal{I}}}
\newcommand{\cJ}{\ensuremath{\mathcal{J}}}
\newcommand{\cK}{\ensuremath{\mathcal{K}}}
\newcommand{\cP}{\ensuremath{\mathcal{P}}}
\newcommand{\cS}{\ensuremath{\mathcal{S}}}

\begin{document}

\begin{center}
{\Large \bf 
Infinite convolutions of probability measures on Polish semigroups 
}
\end{center}
\begin{center}
Kouji \textsc{Yano}\footnote{
Graduate School of Science, Kyoto University.}\footnote{
The research of Kouji Yano was supported by 
JSPS KAKENHI grant no.'s JP19H01791 and JP19K21834 
and by JSPS Open Partnership Joint Research Projects grant no. JPJSBP120209921. 
This research was supported by RIMS and by ISM.}
\end{center}
\begin{center}
{\small \today}
\end{center}

\begin{abstract}
This expository paper is intended for a short self-contained introduction 
to the theory of infinite convolutions of probability measures on Polish semigroups. 
We give the proofs of the Rees decomposition theorem of completely simple semigroups, 
the Ellis--\.{Z}elazko theorem, 
the convolution factorization theorem of convolution idempotents, 
and the convolution factorization theorem of cluster points of infinite convolutions. 
\end{abstract}


\section{Introduction} \label{sec: intro}

As a natural generalization of random walks on an integer lattice, 
the theory of infinite convolutions of probability measures on topological semigroups 
has been extensively studied and widely applied to various problems. 
For this theory, there are celebrated textbooks 
Rosenblatt \cite{Ros}, 
Mukherjea--Tserpes \cite{MT} and H{\"o}gn{\"a}s--Mukherjea \cite{HM}, 
which include a lot of applications of the theory; 
see also Mukherjea's lecture notes \cite{Muk} for applications to random matrices, 
and Ito--Sera--Yano \cite{ISY} for applications to the problem of resolution of $ \sigma $-fields. 

The aim of this paper is to help the reader 
to gain the basic knowledge of this thoery conveniently. 
We mainly follow \cite{HM} and we make some modifications on the proofs. 
For a potential application, 
we develop the theory for topological semigroups with a Polish topology, 
while the textbooks \cite{Ros,MT,HM} deal 
with those with a locally compact Hausdorff second countable topology. 

The goal of this paper is 
the convolution factorization theorem of cluster points of infinite convolutions, 
which will be stated as Theorem \ref{thm: inf conv}. 
The key to the proof is the convolution factorization theorem of convolution idempotents, 
which will be stated as Theorem \ref{thm: conv idem}, 
and the study of probability measures with convolution invariance, 
which will be stated as Proposition \ref{prop: conv inv}. 
Theorems \ref{thm: conv idem} and \ref{thm: inf conv} 
are based on the product decomposition theorem for completely simple semigroups, 
which will be called the \emph{Rees decomposition} and stated as Theorem \ref{thm: Rees decomp}. 
To show that the algebraic decomposition is compatible with a Polish topology, 
we need Ellis--\.{Z}elazko theorem, 
which will be stated as Theorem \ref{thm: Ellis}. 

The Ellis theorem \cite{MR83681}(1957) asserts that 
an algebraic group where the product mapping is separately continuous 
is a topological group, where the topology is locally compact Hausdorff second countable. 
It was extended by \.{Z}elazko \cite{MR125901}(1960) for a completely metrizable topology. 

The study of infinite convolutions on compact groups 
dates back to Kawada--It\^{o} \cite{MR3462}(1940), 
which was generalized by 
Urbanik \cite{MR92921}(1957), 
Kloss \cite{MR123348}(1959), 
and Stromberg \cite{MR114874}(1960) 
and for locally compact groups by 
Tortrat \cite{MR175168}(1964) and Csisz\'{a}r \cite{MR205306}(1966). 
The convolution invariance Proposition \ref{prop: conv inv} 
is due to Mukherjea \cite{MR293687}(1972), 
which originates from the Choquet--Deny equation \cite{MR119041}(1960); 
for later studies, see 
\cite{MR679395, %Wo\'{s} 1982
MR890374, %Raugi 1984
MR826359, %Derriennic 1985
MR876261, %Davies--Shanbhag 1987
MR1010824, %Rao--Shanbhag 1989
MR1083341, %Lau--Zeng 1990
MR1046339}. %Sz\'{e}kely--Zeng 1990
Theorem \ref{thm: conv idem} for convolution idempotents 
is due to Mukherjea--Tserpes \cite{MR296207} (1971); 
for ealier studies, see 
Collins \cite{MR136679}(1962), 
Pym \cite{MR148793}(1962), 
Heble--Rosenblatt \cite{MR169971}(1963), 
Schwarz \cite{MR169970}(1964), 
Choy \cite{MR269770}(1970), 
Duncan \cite{MR276400}(1970), 
and Sun--Tserpes \cite{MR272935}(1970); 
see also \cite{MR1039473}. 
Theorem \ref{thm: inf conv} for cluster points of infinite convolutions 
is due to 
Rosenblatt \cite{MR0118773}(1960) in the compact case 
and to Mukherjea \cite{MR556681}(1979) in the locally compact case; 
for studies earlier than \cite{MR556681}, see 
Glicksberg \cite{MR108690}(1959), 
Collins \cite{MR137789}(1962), 
Schwarz \cite{MR169969}(1964), 
Rosenblatt \cite{MR185636}(1965), 
Lin \cite{MR199306}(1966), 
Mukherjea \cite{MR423458}(1977), 
and Mukherjea--Sun \cite{MR516743}(1978); 
for related papers, see 
\cite{MR543577, %Nakassis 1979
MR841094, %Rosenblatt 1986
MR1135263, %Lo--Mukherjea 1991
MR2050898}. %Budzban--Mukherjea 2004

This paper is organized as follows. 
In Section \ref{sec: algebra} we review the theory of algebraic semigroups. 
In Section \ref{sec: top} we study the theory of Polish semigroups, 
where the Ellis--\.{Z}elazko theorem is proved and utilized. 
Section \ref{sec: conv} is devoted to 
the convolution factorization theorems 
of convolution idempotents and 
of cluster points of infinite convolutions.

\subsection*{Acknowledgements}
The author would like to express our gratitude to Takao Hirayama 
for having a hard time to learn the theory together as beginners. 
He also thanks Yu Ito and Toru Sera for fruitful discussions.

\section{Algebraic semigroup} \label{sec: algebra}

We say that a non-empty set $ S $ is a \emph{semigroup} if 
it is endowed with multiplication 
\begin{align}
S \times S \ni (a,b) \mapsto ab \in S 
\label{}
\end{align}
which is associative, i.e., 
\begin{align}
(ab)c = a(bc) 
, \quad a,b,c \in S. 
\label{}
\end{align}
For two subsets $ A $ and $ B $ of $ S $, 
we denote their product by 
\begin{align}
AB = \{ ab : a \in A , \ b \in B \} . 
\label{}
\end{align}
We write $ A^1 = A $ and $ A^n = A^{n-1} A $ for $ n \ge 2 $. 
We sometimes identify an element $ a \in S $ with the singleton $ \{ a \} $; 
for instance, $ a S = \{ a \} S = \{ ab : b \in S \} $. 
An element $ e \in S $ is called \emph{identity} if 
\begin{align}
xe = ex = x 
, \quad x \in S . 
\label{}
\end{align}
It is obvious that identity is unique if it exists. 
For a semigroup $ S $ with identity $ e $, 
we say that $ y \in S $ is called the \emph{inverse} of $ x \in S $ 
if $ xy = yx = e $. 
It is obvious that the inverse of an element $ x \in S $ is unique if it exists. 
A \emph{group} is a semigroup $ S $ with identity 
such that all elements have their inverses.

\subsection{Left and right simplicity}

Let $ S $ be a semigroup. 
A non-empty subset $ I $ is called a \emph{left ideal} [\emph{right ideal}] (of $ S $) 
if $ S I \subset I $ [$ I S \subset I $]. 
If $ S $ contains no proper left ideal [right ideal], 
then it is called \emph{left simple} [\emph{right simple}]. 
A non-empty subset $ I $ is called a \emph{ideal} if 
it is both a left and a right ideal, i.e., $ S I \cup I S \subset I $. 
If $ S $ contains no proper ideal, then it is called \emph{simple}. 
Note that 
both left and right simple implies simple, but the converse is not necessarily true.

\begin{Lem} \label{lem: minimal left ideal} 
For a subsemigroup $ S $ of a semigroup $ S_0 $, the following are equivalent: 
\begin{enumerate}

\item 
$ S $ is a minimal left ideal of $ S_0 $. 

\item 
$ S = S_0 a $ for all $ a \in S $. 

\end{enumerate}
\end{Lem}

\Proof{
Suppose $ S $ is a minimal left ideal. 
Since $ S_0 a $ for $ a \in S $ is a left ideal of $ S_0 $ 
contained in $ S $, 
we have $ S = S_0 a $ by minimality. 

Suppose $ S = S_0 a $ for all $ a \in S $. 
Let $ I $ be a left ideal of $ S_0 $ such that $ I \subset S $. 
For any $ a \in I $, we have 
$ S = S_0 a \subset S_0 I \subset I $, which shows that 
$ S $ is a minimal left ideal of $ S_0 $. 
}

\begin{Lem} \label{lem: minimal ideal} 
For a subsemigroup $ S $ of a semigroup $ S_0 $, the following are equivalent: 
\begin{enumerate}

\item 
$ S $ is a minimal ideal of $ S_0 $. 

\item 
$ S = S_0 a S_0 $ for all $ a \in S $. 

\end{enumerate}
\end{Lem}

The proof of Lemma \ref{lem: minimal ideal} is almost the same 
as that of Lemma \ref{lem: minimal left ideal}, and so we omit it.

\begin{Lem} \label{lem: left simple}
For a semigroup $ S $, the following are equivalent: 
\begin{enumerate}

\item 
$ S $ is left simple. 

\item 
For any semigroup $ S_0 $ of which $ S $ is a left ideal, 
$ S $ is a minimal left ideal of $ S_0 $. 

\item 
$ S $ is a minimal left ideal of $ S $ itself 
(if and only if $ S = Sa $ for all $ a \in S $ by Lemma \ref{lem: minimal left ideal}). 

\item 
There exists a semigroup $ S_0 $ 
such that $ S $ is a minimal left ideal of $ S_0 $. 

\item 
For any $ a,b \in S $, 
the equation $ xa=b $ has at least one solution $ x \in S $. 

\end{enumerate}
\end{Lem}

\Proof{
[(i) $ \Rightarrow $ (ii)] 
Suppose that $ S $ is a left ideal of a semigroup $ S_0 $ 
and let $ I $ be a left iedal of $ S_0 $ such that $ I \subset S $. 
Then $ S I \subset S_0 I \subset I $, and so $ I $ is a left ideal of $ S $. 
Since $ S $ is left simple, we have $ I = S $, 
which shows that $ S $ is a minimal left ideal of $ S_0 $. 

[(ii) $ \Rightarrow $ (iii) $ \Rightarrow $ (iv)] 
These are obvious. 

[(iv) $ \Rightarrow $ (i)] 
Suppose that $ S $ is a minimal left ideal of $ S_0 $ 
and let $ I $ be a left ideal of $ S $. 
Since $ S_0SI \subset SI \subset I \subset S $, 
we see that $ SI $ is a left ideal of $ S_0 $ with $ SI \subset S $. 
Hence $ SI = S $ by minimality. 
Since $ I \subset S = SI \subset I $, we have $ I = S $, 
which implies that $ S $ is left simple. 

[(iii) $ \Rightarrow $ (v)] 
This is obvious by $ S \subset Sa $. 

[(v) $ \Rightarrow $ (iii)] 
Let $ a \in S $. Then we have $ S \subset Sa $ by (v). 
Since $ S $ is a semigroup, we have $ Sa \subset S $. 
Hence we have $ S = Sa $. 
}

For the simplicity, we have the following. 

\begin{Lem} \label{lem: simple} 
For a semigroup $ S $, the following are equivalent: 
\begin{enumerate}

\item 
$ S $ is simple. 

\item 
For any semigroup $ S_0 $ of which $ S $ is an ideal, 
$ S $ is a minimal ideal of $ S_0 $. 

\item 
$ S $ is a minimal ideal of $ S $ itself 
(if and only if $ S = SaS $ for all $ a \in S $ by Lemma \ref{lem: minimal ideal}). 

\item 
There exists a semigroup $ S_0 $ 
such that $ S $ is a minimal ideal of $ S_0 $. 

\item 
For any $ a,b \in S $, 
the equation $ xay=b $ has at least one solution $ (x,y) \in S \times S $. 

\end{enumerate}
\end{Lem}

The proof of Lemma \ref{lem: simple} 
is almost the same as that of Lemma \ref{lem: left simple}, 
and so we omit it.

\begin{Prop} \label{prop: left right simple}
A semigroup $ S $ which is both left and right simple is a group. 
\end{Prop}

\Proof{
Let $ a \in S $. By Lemma \ref{lem: left simple}, 
we have $ ea = a $ for some $ e \in S $. 
For any $ x \in S $, we have $ x = ay $ for some $ y \in S $, 
and so we have $ ex = eay = ay = x $. 
Similarly, there exists $ e' \in S $ such that $ xe' = x $ for all $ x \in S $. 
Then we obtain $ e' = ee' = e $, and thus $ e $ is identity of $ S $. 

Let $ x \in S $. By Lemma \ref{lem: left simple}, 
we have $ xy = e $ and $ y'x = e $ for some $ y,y' \in S $. 
Since $ y' = y'e = y'xy = ey = y $, we see that $ y $ is the inverse of $ x $. 
}

\subsection{Left and right groups}

Let $ S $ be a semigroup. 
An element $ e \in S $ is called an \emph{idempotent} if $ e^2 = e $. 
We denote the set of all idempotents of $ S $ by 
\begin{align}
E(S) = \{ e \in S : e^2 = e \} . 
\label{}
\end{align}
Note that, if $ e $ is an idempotent, then any element of $ Se $ 
is invariant under right multiplication by $ e $, 
i.e., $ x \in Se $ implies $ xe = x $. 
A semigroup $ S $ is called \emph{left group} [\emph{right group}] 
if $ S $ is left simple [\emph{right simple}] and contains at least one idempotent. 

A semigroup $ S $ is called \emph{left cancellative} [\emph{right cancellative}] 
if, for any $ a,x,y \in S $ with $ ax=ay $ [$ xa=ya $], we have $ x=y $. 
An element $ e \in S $ is called a \emph{left identity} [\emph{right identity}] 
if $ ex=x $ [$ xe=x $] for all $ x \in S $. 

\begin{Lem} \label{lem: right identity}
Let $ S $ be a semigroup. 
If $ S $ is either right cancellative or left simple, 
then any idempotent of $ S $ is a right identity. 
\end{Lem}

\Proof{
Suppose $ S $ is right cancellative and let $ e \in E(S) $. 
Then $ xee = xe $ implies $ xe = x $. 

Suppose $ S $ is left simple and let $ e \in E(S) $. 
By Lemma \ref{lem: left simple}, 
we have $ S = Se $, which yields that 
$ xe=x $ for all $ x \in S $. 
}

\begin{Prop}
For a semigroup $ S $, the following are equivalent: 
\begin{enumerate}

\item 
$ S $ is a left group. 

\item 
$ S $ is left simple and right cancellative. 

\item 
For any $ a,b \in S $, the equation $ xa=b $ 
has a unique solution $ x \in S $. 

\end{enumerate}
\end{Prop}

\Proof{
[(i) $ \Rightarrow $ (ii)] 
Let $ e \in E(S) $ be fixed. 
By Lemma \ref{lem: right identity}, we see that 
$ e $ is a right identity. 

Suppose $ xa=ya $. By Lemma \ref{lem: left simple}, 
we have $ ba = e $ for some $ b \in S $. 
We then have $ abab = aeb = ab $, so that $ ab \in E(S) $ and $ ab $ is a right identity. 
We then obtain $ x = xab = yab = y $. 

[(ii) $ \Rightarrow $ (iii)] 
Existence follows from left simplicity and Lemma \ref{lem: left simple}. 
Uniqueness follows from right cancellativity. 

[(iii) $ \Rightarrow $ (i)] 
By (iii), we have $ S=Sa $ for all $ a \in S $, 
which shows by Lemma \ref{lem: left simple} that $ S $ is left simple. 

Let $ a \in S $ and take $ e \in S $ such that $ ea = a $ by (iii). 
Then we have $ e^2a = ea = a $, which leads to $ e^2 = e $ by right cancellativity. 
}

\subsection{Rees decomposition}

Let $ S $ be a semigroup. 
An idempotent $ e \in E(S) $ is called \emph{primitive} if 
\begin{align}
\text{$ ex=xe=x \in E(S) $ implies $ x=e $}. 
\label{}
\end{align}
We say that $ S $ is \emph{completely simple} 
if $ S $ is simple and contains a primitive idempotent. 

\begin{Thm}[Rees decomposition] \label{thm: Rees decomp} 
Let $ S $ be a completely simple semigroup 
and let $ e $ be a primitive idempotent of $ S $. 
Set 
\begin{align}
L := E(Se) 
, \quad 
G := eSe 
, \quad 
R := E(eS) . 
\label{}
\end{align}
Then the following assertions hold: 
\begin{enumerate}

\item 
$ LG = Se $ is a left group and $ GR = eS $ is a right group. 

\item 
$ RL \subset G $ and 
$ eL = Re = \{ e \} $. 

\item 
$ G = Se \cap eS $ is a group where $ e $ is its identity. 

\item 
$ S = LGR $ (This factorization will be called the \emph{Rees decomposition} of $ S $ at $ e $, 
and $ G $ will be called the \emph{group factor} at $ e $). 

\item 
The product mapping 
\begin{align}
\psi : L \times G \times R \ni (x,g,y) \mapsto (xgy) \in LGR 
\label{}
\end{align}
is bijective with its inverse given as 
\begin{align}
\psi^{-1} : LGR \ni z \mapsto (ze(eze)^{-1},eze,(eze)^{-1}ez) \in L \times G \times R . 
\label{}
\end{align}
\end{enumerate}
\end{Thm}

\Proof{
(i) 
It is obvious that $ Se $ is a left ideal of $ S $. 
Let $ I $ be a left ideal of $ S $ such that $ I \subset Se $. 
Let $ a \in I $. 
Note that $ ae=a $ since $ a \in Se $. 
By simplicity of $ S $ and Lemma \ref{lem: simple}, we have 
$ uav=e $ for some $ u,v \in S $. 
Set $ r = eu $ and $ s = eve $. We then have 
\begin{align}
ras = eu(ae)ve = euave = e 
, \quad 
er = r 
, \quad 
es = se = s . 
\label{}
\end{align}
If we set $ t = sra $, then $ et=te=t $ and 
\begin{align}
t^2 = s(ras)ra = sera = sra = t , 
\label{}
\end{align}
which yields $ t=e $ by primitivity. 
Since $ e = t = sra \in srI \subset I $, 
we have $ Se \subset SI \subset I $, which shows $ I = Se $ 
and that $ Se $ is a minimal left ideal of $ S $. 
By Lemma \ref{lem: left simple}, we see that $ Se $ is left simple. 
Since $ Se $ contains an idempotent $ e $, 
we see that $ Se $ is a left group. 
By a similar argument we see that $ eS $ is a right group. 

Let us show $ LG = Se $. 
It is obvious that $ LG \subset Se $. 
Let $ a \in Se $. Set $ g := ea \in eSe = G $ and set $ b = a g^{-1} \in Se $. 
Since $ g^{-1} = g^{-1} e $, we have 
\begin{align}
b^2 = a g^{-1} a g^{-1} = a g^{-1} (ea) g^{-1} = a g^{-1} = b . 
\label{}
\end{align}
Hence we have $ b \in E(Se) = L $ and $ a = ae = a g^{-1} g = bg \in LG $. 
We now have $ LG = Se $. We also have $ GR = eS $ similarly. 

(ii) 
$ RL \subset (eS)(Se) \subset eSe = G $. 

Let $ x \in L = E(Se) $. 
Since $ (ex)^2 = e(xe)x = exx = ex $ 
and $ e(ex) = (ex)e = ex $, 
we have $ ex=e $ by primitivity. 
We thus see that $ eL = \{ e \} $. 
We have $ Re = \{ e \} $ similarly. 

(iii) 
It is obvious that $ G = eSe = eS \cap Se $, 
since $ x \in eS \cap Se $ implies $ x = ex = xe = exe $. 
It is also obvious that $ e $ is identity of $ G $. 
Let $ g \in G $. Since $ G \subset eSe $, we have $ g = ea $ for some $ a \in Se $. 
By the left simplicity of $ Se $ and by Lemma \ref{lem: left simple}, 
we have $ ba = e $ for some $ b \in Se $. 
Since $ (ab)^2 = a(ba)b = aeb = ab $, 
we see by Lemma \ref{lem: right identity} that $ ab $ is right identity. 
Hence $ ab = abe = e $, which shows that $ b $ is the inverse of $ a $. 

(iv) 
$ LGR = LGGR = SeeS = SeS = S $ by Lemma \ref{lem: simple}. 

(v) 
Let $ z = xgy $ with $ (x,g,y) \in L \times G \times R $. 
Since $ x = xx = xex $ and since $ exgye \in eSe = G $, we have 
\begin{align}
x = xe = x(exgye) (exgye)^{-1} = ze (eze)^{-1} . 
\label{}
\end{align}
We have $ y = (eze)^{-1} ez $ similarly. 
Since $ ex = ye = e $ by (ii), 
we obtain 
\begin{align}
g = ege = (ex)g(ye) = eze. 
\label{}
\end{align}
The proof is now complete. 
}

\begin{Cor} \label{cor: Rees minimal left ideal}
Under the same assumptions and notation as Theorem \ref{thm: Rees decomp}, 
it holds that $ \{ Sy = LGy : y \in R \} $ is the family of 
all minimal left ideals of $ S $. 
\end{Cor}

\Proof{
Any minimal left ideal of $ S $ is of the form $ Sz $ for some $ z \in S $. 
We represent $ z = xgy $ and then we obtain $ Sz = LG(Rx)gy = LGy $, since $ RL \subset G $. 

Conversely, for any $ z \in LGy $, we have $ z = xgy $ for some $ (x,g) \in L \times G $, 
so that we have $ LGyz = LG(yx)gy =LGy $, 
which shows by Lemma \ref{lem: minimal left ideal} that 
$ LGy $ is a minimal left ideal. 
}

\begin{Cor}
Under the same assumptions and notation as Theorem \ref{thm: Rees decomp}, 
the following assertions hold: 
\begin{enumerate}

\item 
For $ z = xgy $ with $ (x,g,y) \in L \times G \times R $, 
$ z $ is idemptent if and only if $ g = (yx)^{-1} $. 

\item 
All idempotents of $ S $ are primitive. 

\item 
Let $ e' $ be another idempotent of $ S $ 
and represent it as $ e' = a(ba)^{-1} b $ for $ (a,b) \in L \times R $. 
Let $ S = L'G'R' $ denote the Rees decomposition of $ S $ at $ e' $. Then 
\begin{align}
L'G' = LGb 
, \quad 
G' = aGb 
, \quad 
G'R' = aGR . 
\label{}
\end{align}
\end{enumerate}
\end{Cor}

\Proof{
(i) Suppose $ z^2 = z $. Then $ xgyxgy=xgy $. 
Since $ eL = Re = \{ e \} $, we have $ gyxg = g $, which shows $ g = (yx)^{-1} $. 

Conversely, suppose $ g = (yx)^{-1} $. 
Then $ z^2 = x(gyxg)y = xgy = z $. 

(ii) Let $ e_1,e_2 \in S $ be two idempotents of $ S $ 
and represent them as $ e_i = a_i (b_i a_i)^{-1} b_i $ 
for $ (a_i,b_i) \in L \times R $, $ i=1,2 $. 
Suppose $ e_1 e_2 = e_2 e_1 = e_2 $. 
Then $ a_1 ((b_1 a_1)^{-1} (b_1a_2) (b_2 a_2)^{-1}) b_2 
= a_2 ((b_2 a_2)^{-1} (b_2 a_1) (b_1 a_1)^{-1}) b_1 = a_2 (b_2 a_2)^{-1} b_2 $, 
which shows $ a_1 = a_2 $ and $ b_1 = b_2 $ 
by the injectivity of the product mapping $ \psi $. 
Hence we have $ e_1 = e_2 $, which shows that $ e_1 $ is a primitive idempotent. 

(iii) 
We have $ L'G' = Se' = LG(Ra)(ba)^{-1}b = LGb $ 
and $ G'R' = aGR $ similarly. 
We also have $ G' = e'Se' = a(ba)^{-1}(bL)G(Ra)(ba)^{-1}b = aGb $. 
}

\begin{Cor}
A left group $ S $ is completely simple. 
The Rees decomposition of $ S $ at $ e \in E(S) $ 
is given as $ S = LG $ with $ R = \{ e \} $. 
\end{Cor}

\Proof{
Suppose $ ex = xe = x \in E(S) $. 
By Lemma \ref{lem: left simple}, 
we have $ yx = e $ for some $ y \in S $. 
Hence $ x = ex = yxx = yx = e $, 
which shows that $ e $ is an primitive idempotent. 
Hence $ S $ is completely simple. 
Let $ S = LGR $ denote the Rees decomposition of $ S $ at $ e $. 
Since $ S=Se $ by Lemma \ref{lem: left simple} 
and since $ Re = \{ e \} $, 
we obtain $ S = Se = LGRe = LG $. 
}

For later use we prove the following proposition. 

\begin{Prop} \label{prop: minimal left right ideal}
Suppose that a semigroup $ S $ contains a minimal left ideal $ A $ 
and a minimal right ideal $ B $ as well. 
Then $ BA $ is a group 
and its identity is a primitive idempotent of $ S $. 
If, in addition, $ S $ is simple, then $ S $ is completely simple. 
\end{Prop}

\Proof{
Since $ (BA)(BA) = (BAB)A \subset BA $, we see that $ BA $ is a subsemigroup of $ S $. 
To prove right simplicity of $ BA $, 
let $ I $ be a right ideal of $ BA $. 
Since $ IB $ is a right ideal of $ S $ and $ IB \subset BAB \subset B $, 
we see that $ IB = B $ by minimality. 
Hence $ BA = IBA \subset I $, which shows right simplicity of $ BA $. 
By a similar argument we obtain left simplicity of $ BA $. 
We thus conclude by Proposition \ref{prop: left right simple} that $ BA $ is a group. 

Let $ e $ be the identity of $ BA $ and suppose $ ex=xe=x \in E(S) $. 
Then $ x = xx = exxe \in (BAS)(SBA) \subset BA $. 
Since $ BA $ is a group and since $ x^2 = x $, we have $ x = xx^{-1} = e $, 
which shows that $ e $ is a primitive idempotent of $ S $. 
}

\subsection{Kernel}

A minimal ideal of a semigroup $ S $ will be called a \emph{kernel} of $ S $. 

\begin{Thm} \label{thm: kernel} 
Let $ S $ be a semigroup. Then the following assertions hold: 
\begin{enumerate}

\item 
If $ S $ contains a minimal left ideal, 
then $ S $ contains a unique kernel $ K $, and 
$ SzS = K $ for all $ z \in K $. 

\item 
If $ S $ contains a minimal left ideal and a minimal right ideal as well, 
then the unique kernel of $ S $ is completely simple. 

\item 
If $ S $ contains a completely simple kernel $ K $, 
then it is the unique kernel of $ S $. 
Let $ K = LGR $ denote the Rees decomposition at $ e $. 
Then $ Sz = Kz = LGz $ for all $ z \in K $. 

\end{enumerate}
\end{Thm}

\Proof{
(i) 
Let $ \cA $ denote the family of all minimal left ideals of $ S $ 
and suppose $ \cA $ is not empty. We shall prove that 
$ K := \bigcup \cA $ is a unique kernel of $ S $. 

Let $ z \in K $ and take $ A \in \cA $ such that $ z \in A $. 
Then $ Sz = A $ by Lemma \ref{lem: minimal left ideal}. 
For $ x \in S $, we see that $ Ax \in \cA $; in fact, 
for any left ideal $ I $ of $ S $ such that $ I \subset Ax $, 
we see that $ J = \{ a \in A : ax \in I \} \subset A $ is a left ideal of $ S $, 
so that $ J=A $ by minimality and thus $ I = Ax $. 
Hence $ SzS = AS = \bigcup_{x \in S} Ax \subset \bigcup \cA = K $, 
which shows by Lemma \ref{lem: minimal ideal} that $ K $ is a kernel of $ S $. 

Let $ K' $ be another kernel of $ S $. 
Since $ K \cap K' $ contains $ KK' $ which is not empty, 
we see that $ K \cap K' $ is an ideal contained both in $ K $ and in $ K' $. 
Thus $ K \cap K' = K = K' $ by minimality. 

(ii) 
By (i) and Lemma \ref{lem: minimal ideal}, we see that 
the unique kernel $ K $ of $ S $ is both a minimal left ideal of $ K $ 
and a minimal right ideal of $ K $. 
By Proposition \ref{prop: minimal left right ideal}, 
we see that $ K $ is completely simple. 

(iii) 
Suppose $ K $ is a completely simple kernel of $ S $ with a primitive idempotent $ e $. 
By Theorem \ref{thm: Rees decomp}, 
$ K $ contains a left group $ Ke $. 
By Lemma \ref{lem: left simple}, 
we see that $ Ke $ is a minimal left idal of $ S $. 
Hence by (i) the kernel of $ S $ is unique. 

For $ z \in K $, we represent $ z = xgy \in LGR $. 
Then by (i) $ Sz $ is a minimal left ideal of $ K $ containing $ y $. 
By Corollary \ref{cor: Rees minimal left ideal}, 
we see that $ Sz = LGy = LGz = Kz $. 
}

\section{Topological semigroup} \label{sec: top}

A semigroup $ S $ is called \emph{topological} 
if $ S $ is endowed with a topology such that 
the product mapping $ S \times S \ni (x,y) \mapsto xy \in S $ is jointly continuous. 
A semigroup $ S $ is called \emph{Polish} 
if $ S $ is a topological semigroup 
with respect to a Polish topology, i.e. a separable and completely metrizable topology. 

It is well-known (see, e.g. \cite[Theorem 1.5.3]{Kechris}) that 
locally compact Polish 
is equivalent to 
locally compact Hausdorff with a countable base. 
It is elementary that 
compact Polish is equivalent to compact metrizable. 

For $ a \in S $ and $ A \subset S $, we write 
\begin{align}
a^{-1} A = \{ x \in S : ax \in A \} 
, \quad 
A a^{-1} = \{ x \in S : xa \in A \} . 
\label{}
\end{align}
If $ S $ contains identity $ e $ and $ a \in S $ has its inverse $ a^{-1} \in S $, then 
$ (a^{-1})A = a^{-1} A $; in fact, 
\begin{align}
(a^{-1})A = \{ a^{-1} x \in S : x \in A \} 
= \{ y \in S : ay \in A \} = a^{-1} A . 
\label{}
\end{align}

\begin{Lem} \label{lem: eA closed}
Let $ S $ be a Polish semigroup. Then the following assertions hold: 
\begin{enumerate}

\item 
For $ a \in S $ and for a closed [open, Borel] subset $ A $, 
both $ a^{-1} A $ and $ A a^{-1} $ are also closed [open, Borel]. 

\item 
If $ A $ is a subsemigroup of $ S $, 
then so is its closure $ \bar{A} $. 

\item 
Let $ A $ be a closed subsemigroup of $ S $. 
Then $ E(A) $, $ eA $, $ Ae $ and $ eAe $ are closed for all $ e \in E(A) $. 

\item 
For two compact subsets $ K $ and $ K' $, 
the product $ KK' $ is also compact. 

\end{enumerate}
\end{Lem}

\Proof{
(i) 
If we write $ \psi_a:S \to S $ for the translation $ \psi_a(x) = ax $, 
then $ a^{-1} A = \psi_a^{-1}(A) $. 
Since $ \psi_a $ is continuous, we obtain the desired results. 

(ii) 
Let $ a,b \in \bar{A} $ and take $ \{ a_n \} , \{ b_n \} \subset A $ 
such that $ a_n \to a $ and $ b_n \to b $. 
Then we have $ ab = \lim a_n b_n \in \bar{A} $. 

(iii) 
Let $ \{ e_n \} \subset E(A) $ such that $ e_n \to e \in S $. 
Since $ A $ is closed, we have $ e \in A $. 
Since $ e_n^2 = e_n $ for all $ n $, we have $ e^2 = e $, which shows $ e \in E(A) $. 

Let $ \{ x_n \} \subset eA $ such that $ x_n \to x \in S $. 
Since $ eA \subset A $ and since $ A $ is closed, 
we have $ x \in A $. Then 
$ ex = \lim ex_n = \lim x_n = x $, which shows $ x = ex \in eA $. 

(iv) 
Let $ \psi:S \times S \to S $ denote the jointly continuous product mapping: 
$ \psi(x,y) = xy $. 
Since $ KK' = \psi(K \times K') $ and $ K \times K' $ is compact, 
we see that $ KK' $ is compact. 
}

\subsection{Topological group}

A group $ S $ is called \emph{topological} 
if $ G $ is a topological semigroup 
and the inverse mapping $ G \ni g \mapsto g^{-1} \in G $ is continuous. 

\begin{Thm}[{Ellis \cite{MR83681} and \.{Z}elazko \cite{MR125901}}] \label{thm: Ellis}
If a group $ G $ is a topological semigroup 
with respect to a completely metrizable topology, 
then it is a topological group. 
\end{Thm}

\Proof{
We borrow the proof from Pfister \cite{MR801345}. 
Let $ e $ denote the identity of $ G $ and 
let $ d $ be a complete metric of $ G $. 

Let $ U_0 $ be a open neighborhood of $ e $. 
By the joint continuity of the product mapping, 
we can construct a sequence $ \{ U_n \}_{n=1}^{\infty } $ of open balls of $ e $ 
such that the radius of $ U_n $ decreases to 0 
and $ \bar{U}_n \bar{U}_n \subset U_{n-1} $ for $ n=1,2,\ldots $, 
where $ \bar{U}_n $ stands for the closure of $ U_n $. 

Let $ \{ x_n \}_{n=1}^{\infty } $ be a subsequence of an arbitrary sequence of $ G $ 
which converges to $ e $. It then suffices to construct a subsequence 
$ \{ n(k) \}_{k=1}^{\infty } $ of $ \{ 1,2,\ldots \} $ such that 
$ x_{n(k)}^{-1} \to e $. 
We write $ y_k := x_{n(1)} \cdots x_{n(k)} $. 

Set $ n(0) = 0 $ and $ y_0 = x_0 = e $. 
If we have $ n(0),n(1),\ldots,n(k-1) $, 
then we can take $ n(k)>n(k-1) $ such that 
$ x_{n(k)} \in U_k $ and $ d(y_k,y_{k-1}) < 2^{-k} $, 
since $ y_{k-1} x_n \to y_{k-1} $ as $ n \to \infty $. 
By completeness of $ d $, we see that $ y_k $ converges to a limit $ y \in G $. 
Let $ n $ be fixed for a while. 
Since $ y U_{n+1} $ is a neighborhood of $ y $, 
we see that $ y_{k-1} \in y U_{n+1} $ for large $ k $. 
For $ j>k $, we have $ U_{j-1} U_j \subset U_j U_j \subset U_{j-1} $, and hence 
\begin{align}
y_k^{-1} y_j = x_{n(k+1)} \cdots x_{n(j-1)} x_{n(j)} 
\in U_{k+1} \cdots U_j \subset U_k , 
\label{}
\end{align}
which implies $ y_k^{-1} y \in \bar{U}_k \subset U_{k-1} $. We now obtain 
\begin{align}
x_{n(k)}^{-1} = (y_{k-1}^{-1} y_k)^{-1} = y_k^{-1} y_{k-1} 
\in y_k^{-1} y U_{n+1} 
\subset U_{k-1} U_{n+1} \subset U_{n+1} U_{n+1} \subset U_n 
\label{}
\end{align}
for large $ k $. Thus we obtain $ x_{n(k)}^{-1} \to e $. 
}

\begin{Cor} \label{cor: Polish semigr csk}
Suppose that a Polish semigroup $ S $ contains a completely simple kernel $ K $. 
Let $ K = LGR $ denote the Rees decomposition of $ K $ at $ e \in E(K) $. 
Then it holds that $ L $, $ G $, $ R $ and $ K $ are closed subsets, 
and that the product mapping 
\begin{align}
\psi : L \times G \times R \ni (x,g,y) \mapsto xgy \in LGR 
\label{}
\end{align}
is homeomorphic. 
\end{Cor}

\Proof{
By Corollary \ref{thm: kernel}, we have $ Ke = Se $, $ eK = eS $ and $ eKe = eSe $. 
By Lemma \ref{lem: eA closed}, we see that 
$ L = E(Ke) $, $ G = eKe $ and $ R = E(eK) $ are all closed. 
By Theorem \ref{thm: Ellis}, we see that $ G $ is a Polish group. 
We now see that the inverse 
\begin{align}
\psi^{-1} : LGR \ni z \mapsto (ze(eze)^{-1},eze,(eze)^{-1}ez) \in L \times G \times R 
\label{}
\end{align}
is continuous. Consequently, we see that $ K $ is closed. 
}

\subsection{Compact semigroup}

\begin{Thm} \label{thm: compact Polish}
A compact Polish semigroup $ S $ contains a compact completely simple kernel. 
\end{Thm}

\Proof{
Let $ \cI $ denote the family of all closed left ideals of $ S $. 
The family $ \cI $ contains $ S $ and is endowed with a partial order by the usual inclusion. 
For any linearly ordered subfamily $ \cJ $ of $ \cI $ has a lower bound in $ \cI $; 
in fact, the intersection $ \bigcap \cJ $ is not empty by compactness of $ S $ 
and is a closed left ideal of $ S $ such that $ \bigcap \cJ \subset J $ for all $ J \in \cJ $. 
Hence, by Zorn's lemma, we see that 
$ \cI $ contains a minimal element, say $ A $. 

Let us prove that $ A $ is a minimal left ideal of $ S $. 
Let $ I $ be a left ideal of $ S $ such that $ I \subset A $. 
For $ a \in I $, we have $ Sa \in \cI $ and $ Sa \subset SI \subset I \subset A $, 
which yields $ Sa = I = A $ by the minimality of $ A $ in $ \cI $. 
This shows that $ A $ is a minimal left ideal of $ S $. 

Similarly we see that $ S $ contains a minimal right ideal. 
By Theorem \ref{thm: kernel}, we see that 
$ S $ contains a completely simple kernel $ K $. 
By Corollary \ref{cor: Polish semigr csk}, we see that $ K $ is a closed subset of $ S $, 
and hence $ K $ is compact. 
}

\begin{Prop} \label{prop: cluster}
Let $ S $ be a Polish semigroup and let $ a \in S $. 
Suppose that 
any subsequence of $ \{ a^n \}_{n=1}^{\infty } $ has a convergent further subsequence. 
Then the set $ C $ of all cluster points of $ \{ a^n \}_{n=1}^{\infty } $ 
is a compact abelian group. 
If we denote the identity of $ C $ by $ e $, then $ C = \bar{\{ e,ae,a^2e,\ldots \}} $. 
\end{Prop}

\Proof{
Let $ C $ denote the set of all cluster points of $ \{ a^n \}_{n=1}^{\infty } $. 
By the assumption, we see that $ C $ is a compact abelian semigroup. 
By Theorem \ref{thm: compact Polish}, we see that 
$ C $ contains a compact completely simple kernel $ K $. 
Since the Rees decomposition of $ K $ is $ LGR = GRL = G $ by commutativity, 
we see that $ K $ is a compact abelian group. 
Let $ e $ denote the identity of $ K $. 
Then, for any $ x \in C $, 
we can find a subsequence $ \{ n(k) \} $ of $ \{ 1,2,\ldots \} $ such that 
$ x = e \lim_{k \to \infty } a^{n(k)} \in eC \subset KC \subset K \subset C $, 
which shows $ K = eC = C $. 
It is now easy to see that $ C = \bar{\{ e,ae,a^2e,\ldots \}} $. 
}

\begin{Rem}
In the settings of Proposition \ref{prop: cluster}, 
suppose that the sequence $ \{ a^n \}_{n=1}^{\infty } $ has multiple points. 
Let $ p $ and $ q $ be the smallest positive integers such that $ a^{q+p} = a^q $. 
Then we have $ \{ a^n : n=1,2,\ldots \} = \{ a,a^2,\ldots,a^{q+p-1} \} $ and 
\begin{align}
K = \{ a^q,a^{q+1},\ldots,a^{q+p-1} \} 
= \{ e,ae,\ldots,a^{p-1}e \} 
\label{}
\end{align}
with $ e = a^{rp} $, where $ r $ is the unique integer such that $ q \le rp \le q+p-1 $. 
\end{Rem}

\section{Convolutions of probability measures on Polish semigroups} \label{sec: conv}

\subsection{Convolutions}

Let $ S $ be a Polish semigroup. 
Let $ \cB(S) $ denote the family of all Borel sets of $ S $ 
and $ \cP(S) $ the family of all probability measures on $ (S,\cB(S)) $. 

For $ \mu,\nu \in \cP(S) $, 
we define the \emph{convolution} $ \mu * \nu \in \cP(S) $ of $ \mu $ and $ \nu $ by 
\begin{align}
\mu * \nu (B) = \iint 1_B(xy) \mu(\d x) \nu(\d y) 
, \quad B \in \cB(S) . 
\label{}
\end{align}
Since $ 1_B(xy) = 1_{By^{-1}}(x) = 1_{x^{-1}B}(y) $, we have 
\begin{align}
\mu * \nu (B) = \int \mu(By^{-1}) \nu(\d y) 
= \int \nu(x^{-1}B) \mu(\d x) 
, \quad B \in \cB(S) . 
\label{}
\end{align}
For $ a \in S $, we write $ \delta_a $ for the Dirac mass at $ a $: 
$ \delta_a(B) = 1_B(a) $. It is obvious that 
\begin{align}
\mu * \delta_x(B) = \mu(B x^{-1}) 
, \quad 
\delta_x * \mu(B) = \mu(x^{-1} B) 
, \quad 
B \in \cB(S) , 
\label{}
\end{align}
which will be called \emph{translations} of $ \mu $.

For $ \mu \in \cP(S) $, we denote its \emph{topological support} by 
\begin{align}
\cS(\mu) = \{ x \in S : \text{$ \mu(U) > 0 $ for all open neighborhood $ U $ of $ x $} \} . 
\label{}
\end{align}
It is obvious that $ \cS(\mu) $ is closed and $ \mu(\cS(\mu)^c) = 0 $. 

\begin{Lem}
For $ \mu,\nu \in \cP(S) $, it holds that 
\begin{align}
\cS(\mu * \nu) = \bar{ \cS(\mu) \cS(\nu) } . 
\label{}
\end{align}
\end{Lem}

\Proof{
Let $ a \in \cS(\mu) $ and $ b \in \cS(\nu) $. 
For any open neighborhood $ U $ of $ ab $, 
the joint continuity of the product mapping allows us to take open neighborhoods 
$ U_1 $ of $ a $ and $ U_2 $ of $ b $ such that $ U_1 U_2 \subset U $, so that 
\begin{align}
\mu * \nu(U) 
\ge \iint 1_{U_1 U_2}(xy) \mu(\d x) \nu(\d y) 
\ge \mu(U_1) \nu(U_2) > 0 , 
\label{}
\end{align}
which yields $ ab \in \cS(\mu * \nu) $ 
and hence $ \bar{ \cS(\mu) \cS(\nu) } \subset \cS(\mu * \nu) $. 

Let $ a \in \bar{ \cS(\mu) \cS(\nu) }^c $. 
Then we can take an open neighborhood $ U $ of $ a $ such that 
$ U \subset \cbra{ \cS(\mu) \cS(\nu) }^c $, so that 
\begin{align}
\mu * \nu(U) 
\le \iint 1_{\cbra{ \cS(\mu) \cS(\nu) }^c}(xy) \mu(\d x) \nu(\d y) 
\le \mu(\cS(\mu)^c) + \nu(\cS(\nu)^c) = 0 , 
\label{}
\end{align}
which shows $ a \in \cS(\mu * \nu)^c $ 
and hence $ \cS(\mu * \nu) \subset \bar{ \cS(\mu) \cS(\nu) } $. 
}

\begin{Prop} \label{prop: mu*nu LR} 
Let $ S $ be a completely simple Polish semigroup. 
Let $ S = LGR $ denote the Rees decomposotion at $ e \in E(S) $. 
For the inverse of the product mapping $ \psi:L \times G \times R \to LGR $, we denote 
\begin{align}
(z^L,z^G,z^R) := \psi^{-1}(z) = (ze(eze)^{-1},eze,(eze)^{-1}ez) \in L \times G \times R 
, \quad z \in LGR . 
\label{}
\end{align}
For $ \mu \in \cP(S) $, we define 
\begin{align}
\mu^L(B) = \mu(z : z^L \in B) 
, \quad 
\mu^G(B) = \mu(z : z^G \in B) 
, \quad 
\mu^R(B) = \mu(z : z^R \in B) 
\label{eq: muL}
\end{align}
for $ B \in \cB(S) $. 
Then, for $ \mu,\nu \in \cP(S) $, it holds that 
\begin{align}
(\mu * \nu)^L = \mu^L 
, \quad 
(\mu * \nu)^R = \nu^R . 
\label{}
\end{align}
\end{Prop}

\Proof{
This is obvious by noting that 
$ (z_1 z_2)^L = z_1^L $ and $ (z_1 z_2)^R = z_2^R $. 
}

We equip $ \cP(S) $ with the topology of weak convergence: 
$ \mu_n \to \mu $ if and only if 
$ \int f \d \mu_n \to \int f \d \mu $ 
for all $ f \in C_b(S) $, 
the class of all bounded continuous functions on $ S $. 
It is well-known (see, e.g. \cite[Theorems 6.2 and 6.5 of Chapter 2]{Par}) that 
$ \cP(S) $ is a Polish space. 

\begin{Prop}
Let $ S $ be a Polish semigroup. 
Then the convolution mapping 
$ \cP(S) \times \cP(S) \ni (\mu,\nu) \mapsto \mu * \nu \in \cP(S) $ 
is jointly continuous. Consequently, $ \cP(S) $ is a Polish semigroup. 
\end{Prop}

\Proof{
Note that, 
if we take independent random variables $ X $ and $ Y $ taking values in $ S $ 
such that $ X \dist \mu $ and $ Y \dist \nu $, then 
 $ \mu * \nu $ coincides with the law of the product $ XY $. 
The desired result now follows from the Skorokhod coupling thoerem 
(see, e.g. \cite[Theorem 4.30]{Kal}), which asserts that 
$ \mu_n \to \mu $ implies that 
we can take random variables $ \{ X_n \},X $ taking values in $ S $ such that 
$ X_n \dist \mu_n $, $ X \dist \mu $ and $ X_n \to X $ a.s. 
}

\subsection{Translation invariance}

Let $ S $ be a Polish semigroup. 
A probability measure $ \mu \in \cP(S) $ is called 
\emph{$ \ell^* $-invariant} [\emph{$ r^* $-invariant}] if 
$ \delta_x * \mu = \mu $ [$ \mu * \delta_x = \mu $] for all $ x \in S $. 

\begin{Thm} \label{thm: ell* r* inv}
Let $ S $ be a Polish semigroup and let $ \mu \in \cP(S) $. 
Suppose that $ \mu $ is both $ \ell^* $-invariant and $ r^* $-invariant. 
Then $ \cS(\mu) $ is a compact Polish group, 
and $ \mu $ coincides with the normalized unimodular Haar measure on $ \cS(\mu) $ 
(see e.g. \cite[Chapter 9]{Coh} for the Haar measure). 
\end{Thm}

\Proof{
Note that 
\begin{align}
\cS(\mu) = \cS(\delta_x * \mu) = \bar{ x \cS(\mu) } 
, \quad x \in S , 
\label{}
\end{align}
which implies that $ \cS(\mu) $ is a left ideal of $ S $. 
Similarly $ \cS(\mu) $ is a right ideal of $ S $, 
and hence $ \cS(\mu) $ is an ideal of $ S $. 

Let us prove that, for any $ x \in \cS(\mu) $, the subsemigroup $ x S $ is left-cancellative. 
Let $ y,a,b \in S $ be such that $ (xy)(xa) = (xy)(xb) $. 
Since $ \cS(\mu) = \cS(\mu * \delta_{xyx}) = \bar{\cS(\mu) xyx} $, 
we can take $ \{ z_n \} \subset \cS(\mu) $ such that $ z_nxyx \to x $, and hence 
\begin{align}
xa = \lim z_n xyxa = \lim z_n xyxb = xb , 
\label{}
\end{align}
which shows that $ x S $ is left-cancellative. 
Similarly $ S x $ is right-cancellative. 

Let $ a,b \in \cS(\mu) $ be fixed. 
We shall prove that the subsemigroup $ D := a \cS(\mu) b $ contains an idempotent. 
Note that 
\begin{align}
\mu(D) = (\delta_a * \mu * \delta_b)(D) = \mu(a^{-1}Db^{-1}) 
\ge \mu(\cS(\mu)) = \mu(S) = 1 , 
\label{}
\end{align}
which shows $ \mu(D) = 1 $. For $ x \in D $, we have 
\begin{align}
\mu(D) \le \mu(x^{-1}(xD)) = (\mu*\delta_x)(xD) = \mu(xD) \le \mu(D) , 
\label{}
\end{align}
which shows $ \mu(xD) = \mu(D) = 1 $. 
We define two mappings $ \theta ,\beta : S \times S \to S \times S $ by 
\begin{align}
\theta(x,y) = (x,xy) 
, \quad 
\beta(x,y) = (y,x) . 
\label{}
\end{align}
Since $ (x,y) \in \theta(D \times D) $ if and only if $ x \in D $ and $ y \in xD $, 
we have 
\begin{align}
(\mu \otimes \mu)(\beta \circ \theta(D \times D)) 
= (\mu \otimes \mu)(\theta(D \times D)) 
= \int_D \mu(xD) \mu(\d x) = \mu(D)^2 = 1 . 
\label{}
\end{align}
This shows that $ \beta \circ \theta(D \times D) \cap \theta(D \times D) $ is not empty, 
so that $ (vw,v) = (x,xy) $ for some $ v,w,x,y \in D $. 
We now have $ x(yw) = vwyw = x(yw)^2 $, 
which implies $ yw=(yw)^2 $ by left-cancellativity of $ D $. 

Let $ e := yw \in E(D) = E(a\cS(\mu)b) $. 
By the left- and right-cancellativity of $ a\cS(\mu)b $ 
and by Lemma \ref{lem: right identity}, 
we see that $ e $ is identity of $ a\cS(\mu)b $. 
By Lemma \ref{lem: eA closed}, we see that 
\begin{align}
\cS(\mu) = \bar{ e a \cS(\mu) b } 
= e \, (\bar{ a \cS(\mu) b }) 
= e \cS(\mu) 
\subset a \cS(\mu) b \cS(\mu) 
\subset a \cS(\mu) \subset \cS(\mu) , 
\label{}
\end{align}
which shows $ a \cS(\mu) = \cS(\mu) $. 
Similarly we have $ \cS(\mu) b = \cS(\mu) $. 
By Lemma \ref{lem: left simple}, Proposition \ref{prop: left right simple} 
and Theorem \ref{thm: Ellis}, 
we see that $ \cS(\mu) $ is a Polish group.

By the $ \ell^* $-invariance, we have $ \mu * \mu = \mu $. 
We now apply \cite[Theorem 3.1 of Chapter 3]{Par} 
to obtain the desired result. 
}

\subsection{Convolution invariance}

\begin{Prop}[{Mukherjea \cite{MR293687}}] \label{prop: conv inv}
Let $ S $ be a Polish semigroup and let $ \mu,\nu \in \cP(S) $. 
Suppose 
\begin{align}
\nu = \mu * \nu = \nu * \mu . 
\label{}
\end{align}
Then, for any $ x \in \cS(\mu) $ and any $ a \in \cS(\nu) $, it holds that 
\begin{align}
\nu * \delta_{xa} = \nu * \delta_a 
, \quad 
\delta_{ax} * \nu = \delta_a * \nu . 
\label{}
\end{align}
\end{Prop}

\Proof{
Let $ a \in \cS(\nu) $, $ f \in C_b(S) $ and $ \eps>0 $ be fixed for a while, and set 
\begin{align}
g(x) = \max \cbra{ \int f \d (\nu * \delta_x) - \int f \d (\nu * \delta_a) - \eps , 0 } 
, \quad x \in S. 
\label{}
\end{align}
It is obvious that $ g \in C_b(S) $, $ g $ is non-negative and $ g(a) = 0 $. 
By $ \nu = \nu * \mu $, we have 
\begin{align}
& \int f \d (\nu * \delta_x) - \int f \d (\nu * \delta_a) - \eps 
\label{} \\
=& \int \cbra{ \int f \d (\nu * \delta_{yx}) - \int f \d (\nu * \delta_a) - \eps } \mu(\d y) 
\le \int g(yx) \mu(\d y) , 
\label{}
\end{align}
so that we have 
\begin{align}
g(x) \le \int g(yx) \mu(\d y) 
, \quad x \in S . 
\label{eq: ineq gx}
\end{align}
In addition, by $ \nu = \mu * \nu $, we have 
\begin{align}
\int \cbra{ g(x) - \int g(yx) \mu(\d y) } \nu(\d x) 
= \int g \d \nu - \int g \d (\mu*\nu) = 0 , 
\label{}
\end{align}
which shows that the equality in \eqref{eq: ineq gx} holds 
for $ \nu $-a.e. $ x \in S $. 
Since $ g $ is continuous, we see that 
the equality in \eqref{eq: ineq gx} holds 
for all $ x \in \cS(\nu) $. 
Since $ a \in \cS(\nu) $ and $ g(a)=0 $, we see, 
again by continuity of $ g $, that 
\begin{align}
g(ya) = 0 , \quad y \in \cS(\mu) . 
\label{}
\end{align}
Since $ \eps>0 $ is arbitrary, we obtain 
\begin{align}
\int f \d (\nu * \delta_{ya}) \le \int f \d (\nu * \delta_a) 
, \quad a \in \cS(\nu) , \ y \in \cS(\mu) . 
\label{}
\end{align}
Since $ \nu = \nu * \mu $, we have 
$ \int \cbra{ \int f \d (\nu * \delta_{ya}) - \int f \d (\nu * \delta_a) } \mu(\d y) = 0 $, 
which implies 
\begin{align}
\int f \d (\nu * \delta_{ya}) = \int f \d (\nu * \delta_a) 
, \quad 
a \in \cS(\nu) , \ y \in \cS(\mu) , \ f \in C_b(S) . 
\label{}
\end{align}
Since $ f \in C_b(S) $ is arbitrary, we obtain $ \nu * \delta_{ya} = \nu * \delta_a $ 
for all $ a \in \cS(\nu) $ and $ y \in \cS(\mu) $. 
We obtain $ \delta_{ay} * \nu = \delta_a * \nu $ similarly. 
}

\subsection{Convolution idempotent}

We denote the $ n $-fold convolution by $ \mu^n $, i.e. 
$ \mu^1 = \mu $ and $ \mu^n = \mu^{n-1} * \mu $ for $ n=2,3,\ldots $. 

\begin{Thm}[{Mukherjea--Tserpes \cite{MR296207}}] \label{thm: conv idem}
Let $ S $ be a Polish semigroup and let $ \mu \in \cP(S) $. 
Suppose that $ \mu^2 = \mu $. 
Then $ \cS(\mu) $ is completely simple and its group factor is compact. 
Let $ \cS(\mu) = LGR $ denote the Rees decomposition at $ e \in E(\cS(\mu)) $. 
Then $ \mu $ admits the convolution factorization 
\begin{align}
\mu = \mu^L * \omega_G * \mu^R , 
\label{eq: muL oG muR}
\end{align}
where $ \mu^L $ and $ \mu^R $ have been introduced in \eqref{eq: muL} 
and $ \omega_G $ stands for the normalized unimodular Haar measure 
on the compact Polish group $ G $. 
\end{Thm}

\begin{Rem}
The convolution factorization \eqref{eq: muL oG muR} is equivalent to the following assertion: 
If we let $ Z $ be a random variable whose law is $ \mu $, then 
\begin{align}
\text{$ Z^L $, $ Z^G $ and $ Z^R $ are independent and the law of $ Z^G $ is $ \omega_G $}. 
\label{}
\end{align}
Here $ (Z^L,Z^G,Z^R) = \psi^{-1}(Z) $ with $ \psi:L \times G \times R \to LGR $ 
denoting the product mapping; see Proposition \ref{prop: mu*nu LR}. 
\end{Rem}

\Proof[Proof of Theorem \ref{thm: conv idem}]{
Since $ \cS(\mu) = \bar{\cS(\mu)\cS(\mu)} $, 
we see that $ \cS(\mu) $ is a closed subsemigroup of $ S $. 
By Proposition \ref{prop: conv inv}, we see that, for any $ a \in \cS(\mu) $, 
\begin{align}
\mu * \delta_{xa} = \mu * \delta_a 
, \quad 
\delta_{ax} * \mu = \delta_a * \mu 
, \quad x \in \cS(\mu) . 
\label{eq: axmu}
\end{align}
Then, for $ a \in \cS(\mu) $, we have 
\begin{align}
\mu * \delta_{ay} = \mu * \delta_a \ (y \in \cS(\mu * \delta_a)) . 
, \quad 
\delta_{za} * \mu = \delta_a * \mu \ (z \in \cS(\delta_a * \mu)) 
\label{eq: yamu}
\end{align}
In fact, for $ y \in \cS(\mu * \delta_a) = \bar{\cS(\mu) a} $, 
we may take $ \{ x_n \} \subset \cS(\mu) $ such that $ x_n a \to y $, so that 
$ \mu * \delta_a = \mu * \delta_{ax_na} \to \mu * \delta_{ay} $. 

Let $ a \in \cS(\mu) $ be fixed and set $ \nu = \delta_a * \mu * \delta_a $. 
Then $ \cS(\nu) = \bar{a \cS(\mu) a} $ is a closed subsemigroup of $ S $. 
For any $ y \in \cS(\nu) = \bar{a \cS(\mu) a} $, 
we may take $ \{ x_n \} \subset \cS(\mu) $ such that $ a x_n a \to y $, 
so that, using \eqref{eq: axmu}, we have 
\begin{align}
\nu = \delta_a * \mu * \delta_a = \delta_{a x_n a^2} * \mu * \delta_a 
= \delta_{a x_n a} * \nu \to \delta_y * \nu , 
\label{}
\end{align}
which shows that $ \nu|_{\cS(\nu)} $ is $ \ell^* $-invariant. 
We see similarly that $ \nu|_{\cS(\nu)} $ is $ r^* $-invariant. 
We may now apply Theorem \ref{thm: ell* r* inv} 
to see that $ \cS(\nu) = \bar{a \cS(\mu) a} $ is a compact Polish group. 
Its identity is an idempotent of $ \cS(\mu) $. 

Let $ e \in E(\cS(\mu)) $. 
By the above argument with $ a=e $, 
we see that $ G := e \cS(\mu) e $ is a compact Polish group 
(note that $ e \cS(\mu) e $ is closed by Lemma \ref{lem: eA closed}). 
Set $ A := \cS(\mu) e $. 
For $ y \in A $, using \eqref{eq: yamu}, we have 
\begin{align}
\bar{Ay} = \bar{ \cS(\mu) ey }  = \cS(\mu * \delta_{ey}) 
= \cS(\mu * \delta_e ) = \cS(\mu) e = A . 
\label{}
\end{align}
Since $ Ay \cap e \cS(\mu) e $ is a left ideal of the group $ e \cS(\mu) e $, 
we see that $ Ay \cap e \cS(\mu) e = e \cS(\mu) e $, i.e. $ e \cS(\mu) e \subset Ay $, 
which shows $ e \in Ay $. Hence 
\begin{align}
A = Ae \subset AAy \subset Ay \subset \bar{Ay} = A , 
\label{}
\end{align}
which yields $ Ay = A $ for all $ y \in A $. 
By Lemma \ref{lem: left simple}, we see that $ A $ is a left group. 
We see similarly that $ B := e \cS(\mu) $ is a right group. 
By Theorem \ref{thm: kernel}, we see that $ \cS(\mu) $ contains a completely simple kernel $ K $, 
which is closed by Corollary \ref{cor: Polish semigr csk}. 

By \eqref{eq: axmu}, we have 
\begin{align}
\mu * \delta_e * \mu 
= \int (\mu * \delta_e * \delta_a) \mu(\d a) 
= \int (\mu * \delta_a) \mu(\d a) = \mu * \mu = \mu . 
\label{}
\end{align}
By Lemma \ref{lem: simple}, we have $ K = \cS(\mu) e \cS(\mu) $, and hence we obtain 
\begin{align}
K = \bar{K} = \bar{ \cS(\mu) e \cS(\mu) } 
= \cS(\mu * \delta_e * \mu) = \cS(\mu) , 
\label{}
\end{align}
which shows that $ \cS(\mu) $ is completely simple. 

By \eqref{eq: yamu}, we see that $ \mu * \delta_e $ is $ r^* $-invariant 
on $ A = \cS(\mu) e = LG $, so that 
$ \mu * \delta_e = \mu * \delta_e * \omega_G $. Hence, for any $ B \in \cB(\cS(\mu)) $, 
\begin{align}
\mu(B) 
=& (\mu * \delta_e * \mu)(B) 
= (\mu * \delta_e * \omega_G * \mu)(B) 
\label{} \\
=& \int \mu(\d z_1) \int \mu(\d z_2) \int \omega_G(\d g) 1_B(z_1 e g z_2) 
\label{} \\
=& \int \mu(\d z_1) \int \mu(\d z_2) \int \omega_G(\d g) 1_B(z_1^L g z^R_2) 
= (\mu^L * \omega_G * \mu^R)(B) , 
\label{}
\end{align}
which completes the proof. 
}

The following proposition is a converse to Theorem \ref{thm: conv idem}. 

\begin{Prop}
Let $ S $ be a Polish semigroup and let $ \mu_1,\mu_2 \in \cP(S) $. 
Let $ G $ be a compact Polish subgroup of $ S $ 
and suppose that $ \cS(\mu_2 * \mu_1) \subset G $. 
Then $ \mu := \mu_1 * \omega_G * \mu_2 $ satisfies $ \mu^2 = \mu $. 
\end{Prop}

\Proof{
For any $ B \in \cB(S) $, we have 
\begin{align}
& \mu^2(B) = (\mu_1 * \omega_G * \mu_2 * \mu_1 * \omega_G * \mu_2) (B) 
\label{} \\
=& \int \mu_1(\d z_1) \int \omega_G(\d g_1) \int (\mu_2 * \mu_1)(\d g_2) \int \omega_G(\d g_3) \int \mu_2(\d z_2) 1_B(z_1 g_1 g_2 g_3 z_2) 
\label{} \\
=& \int \mu_1(\d z_1) \int \omega_G(\d g_1) \int \mu_2(\d z_2) 1_B(z_1 g_1 z_2) 
= (\mu_1 * \omega_G * \mu_2) (B) = \mu(B) , 
\label{}
\end{align}
which completes the proof. 
}

\subsection{Infinite convolutions}

\begin{Thm}[{Rosenblatt \cite{MR0118773} and Mukherjea \cite{MR556681}}] \label{thm: inf conv}
Let $ S_0 $ be a Polish semigroup and let $ \mu \in \cP(S_0) $. 
Suppose that the sequence $ \{ \mu^n \}_{n=1}^{\infty } $ is tight. 
Let $ S $ denote the closure of the semigroup generated by $ \cS(\mu) $, i.e. 
\begin{align}
S := \bar{ \bigcup_{n=1}^{\infty } \cS(\mu)^n } . 
\label{}
\end{align}
Then the following assertions hold: 
\begin{enumerate}

\item 
There exists $ \nu \in \cP(S) $ such that 
$ \nu^2 = \nu $, $ \mu * \nu = \nu * \mu = \nu $ and 
\begin{align}
\mu_n := \frac{1}{n} \sum_{k=1}^n \mu^k \tend{}{n \to \infty } \nu . 
\label{}
\end{align}

\item 
The family $ \cK $ of cluster points of $ \{ \mu^n : n=1,2,\ldots \} $ 
is a compact abelian group such that 
\begin{align}
\cS(\nu) = \bar{ \bigcup_{\lambda \in \cK} \cS(\lambda) } . 
\label{}
\end{align}

\item 
Let $ \eta $ denote the identity of $ \cK $. 
Then $ \cS(\eta) $ is a completely simple semigroup. 
Let $ \cS(\eta) = LHR $ denote the Rees decomposition at $ e \in E(\cS(\eta)) $. 
Then $ H $ is a compact group 
and $ \eta $ admits the convolution factorization 
\begin{align}
\eta = \eta_L * \omega_H * \eta_R . 
\label{eq: conv factr eta}
\end{align}

\item 
$ \cS(\nu) $ is a completely simple kernel of $ S $ 
containing the idempotent $ e $. 
The Rees decomposition of $ \cS(\nu) $ at $ e $ is of the form $ \cS(\nu) = LGR $, 
where $ G $ is a compact group containing $ H $, 
and $ \nu $ admits the convolution factorization 
\begin{align}
\nu = \eta_L * \omega_G * \eta_R . 
\label{}
\end{align}

\item 
For $ g \in G $, we write $ \omega_{gH} := \delta_g * \omega_H $. 
It holds that $ H $ is a closed normal subgroup of $ G $ 
and that there exists a Polish group isomorphism $ F : \cK \to G/H $ such that 
\begin{align}
\lambda = \eta^L * \omega_{F(\lambda)} * \eta^R , 
\label{}
\end{align}
Consequently, there exists $ \gamma \in G $ such that 
$ \mu^k * \eta $ admits the convolution factorization 
\begin{align}
\mu^k * \eta = \eta^L * \omega_{\gamma^k H} * \eta^R 
, \quad k=1,2,\ldots, 
\label{eq: conv factr muketa}
\end{align}
and furthermore, $ \cK $ and $ G/H $ may be represented as 
\begin{align}
\cK = \bar{ \{ \eta,\mu*\eta,\mu^2*\eta,\ldots \} } 
, \quad 
G/H = \bar{ \{ H,\gamma H,\gamma^2 H,\ldots \} } . 
\label{eq: cK G repre}
\end{align}
\end{enumerate}
\end{Thm}

\begin{Rem}
If the order of the group $ \cK $ or $ G/H $ is finite, say $ p $, then 
\begin{align}
\cK = \{ \eta,\mu*\eta,\ldots,\mu^{p-1}*\eta \} 
, \quad 
G/H = \{ H,\gamma H,\ldots,\gamma^{p-1} H \} 
\label{}
\end{align}
with $ \gamma^p \in H $. 
It is now obvious that $ \lim_{n \to \infty } \mu^n $ converges if and only if $ p=1 $. 
\end{Rem}

\Proof[Proof of Theorem \ref{thm: inf conv}]{
(i) 
Let $ \| \cdot \| $ denote the total variation norm. 
For $ j=1,2,\ldots $, we have 
\begin{align}
\norm{ \mu_n - \mu^j * \mu_n } 
\le \frac{1}{n} \norm{ \sum_{k=1}^n \mu^k - \sum_{k=1}^n \mu^{k+j} } 
= \frac{1}{n} \norm{ \sum_{k=1}^j \mu^k - \sum_{k=n+1}^{n+j} \mu^k } 
\le \frac{2j}{n} \tend{}{n \to \infty } 0. 
\label{eq: mun-mujmun}
\end{align}
Since $ \{ \mu^n \} $ is tight, we see that $ \{ \mu_n \} $ is also tight. 
Let $ \nu_1,\nu_2 $ be cluster points of $ \{ \mu_n \} $. 
For $ i=1,2 $, we see by \eqref{eq: mun-mujmun} that $ \mu^j * \nu_i = \nu_i * \mu^j = \nu_i $ 
for $ j=1,2,\ldots $, so that 
$ \mu_n * \nu_i = \nu_i * \mu_n = \nu_i $ for $ n=1,2,\ldots $, 
which implies $ \nu_1 = \nu_1 * \nu_2 = \nu_2 * \nu_1 = \nu_2 $. 
Hence we see that $ \{ \mu_n \} $ converges to some $ \nu \in \cP(S_0) $ 
and we have $ \nu^2 = \nu $ and $ \mu * \nu = \nu * \mu = \nu $. 
We may apply Theorem \ref{thm: conv idem} to see that 
$ \cS(\nu) $ is a completely simple semigroup 
and its group factor is compact. 

(ii) 
Let us prove that $ \cS(\nu) $ 
and $ \cS(\cK) := \bar{\bigcup_{\lambda \in \cK} \cS(\lambda) } $ 
are ideals of $ S $. 
Let $ a \in S $, $ x \in \cS(\nu) $ and $ y \in \cS(\cK) $. 
Then we may take $ \{ a_n \} \subset \cS(\mu)^{m(n)} \subset \cS(\mu^{m(n)}) $ 
and $ \{ y_n \} \subset \cS(\lambda_n) $ such that 
$ a_n \to a $ and $ y_n \to y $. 
Since 
\begin{align}
a_n x \in& \cS(\mu^{m(n)}) \cS(\nu) \subset \cS(\mu^{m(n)} * \nu) = \cS(\nu) 
, \label{} \\
a_n y_n \in& \cS(\mu^{m(n)}) \cS(\lambda_n) \subset \cS(\mu^{m(n)} * \lambda_n) \subset \cS(\cK) , 
\label{}
\end{align}
we obtain $ ax = \lim a_n x \in \cS(\nu) $ and $ ay = \lim a_n y_n \in \cS(\cK) $, 
which shows that $ \cS(\nu) $ and $ \cS(\cK) $ are both left ideals of $ S $. 
Similarly we see that they are also right ideals of $ S $. 

Let $ U $ be an open subset containing $ \cS(\nu) $. 
We shall prove that $ \mu^n(U) \to 1 $. 
Let $ \eps>0 $. By tightness, we may take a compact subset $ K_1 $ such that 
$ \inf_n \mu^n(K_1) > 1-\eps $. 
We may take a compact subset $ K_2 \subset \cS(\nu) $ such that 
$ \nu(K_2) > 1-\eps $. 
Since $ K_1 K_2 \subset S \cS(\nu) \subset \cS(\nu) \subset U $, 
we have $ K_1 \times K_2 \subset \tilde{U} := \{ (x,y) \in S_0 \times S_0 : xy \in U \} $. 
By the Wallace theorem (see, e.g., \cite[Theorem 12 of Chapter 5]{Kel}), 
we may take open subsets $ V_1 $ and $ V_2 $ 
such that $ K_1 \subset V_1 $, $ K_2 \subset V_2 $ and $ V_1 \times V_2 \subset \tilde{U} $, 
which implies $ V_1 V_2 \subset U $. 
Since $ \mu_n \to \nu $, we have $ \liminf_n \mu_n(V_2) \ge \nu(V_2) \ge \nu(K_2) > 1-\eps $. 
We may then take some $ n_0 $ such that $ \mu^{n_0}(V) > 1-\eps $. 
We now have 
\begin{align}
\mu^{n+n_0}(U) = \iint 1_U(xy) \mu^n(\d x) \mu^{n_0}(\d y) 
\ge \mu^n(V_1) \mu^{n_0}(V_2) > (1-\eps)^2, 
\label{}
\end{align}
which leads to $ \mu^n(U) \to 1 $.

By the tightness assumption, we may apply Proposition \ref{prop: cluster} 
to see that $ \cK $ is a compact abelian group. 
Let $ \lambda \in \cK $ and let $ x \in S(\lambda) $. 
Suppose that $ x \notin \cS(\nu) $. 
We could then take disjoint open sets $ U $ and $ V $ 
such that $ \cS(\nu) \subset U $ and $ x \in V $. 
If we let $ \delta := \lambda(V)/2 > 0 $, then 
$ \mu^n(V) > \delta $ for infinitely many $ n $, 
and then $ \liminf_n \mu^n(U) \le \liminf_n \mu^n(V^c) \le 1-\delta $, 
which would contradict $ \mu^n(U) \to 1 $. 
Hence we obtain $ \cS(\cK) \subset \cS(\nu) $. 
Since $ \cS(\nu) $ is a minimal ideal of $ S $ by Lemma \ref{lem: simple} 
and since $ \cS(\cK) $ is an ideal of $ S $, 
we see that $ \cS(\cK) = \cS(\nu) $. 

(iii) 
By Theorem \ref{thm: conv idem}, we see that 
$ \cS(\eta) $ is a completely simple semigroup. 
Let $ \cS(\eta) = LHR $ denote the Rees decomposition at $ e \in E(\cS(\eta)) $ 
(hence $ RL \subset H $). 
Then the group factor $ H $ is compact 
and $ \eta $ admits the convolution factorization \eqref{eq: conv factr eta}. 

(iv) 
We have already seen in (i) that $ \cS(\nu) $ is a completely simple kernel of $ S $. 
Since $ \cS(\eta) \subset \cS(\cK) = \cS(\nu) $, we have $ e \in E(\cS(\nu)) $. 
Let $ \cS(\nu) = L'GR' $ denote the Rees decomposition at $ e $. 
As a consequence of Theorem \ref{thm: conv idem}, we see that 
$ \nu $ admits the convolution factorization $ \nu = \eta^{L'} * \omega_G * \eta^{R'} $. 
Since $ \cS(\eta) \subset \cS(\nu) $ and $ L = E(\cS(\eta)e)) $ etc., we see that 
$ L \subset L' $, $ H \subset G $ and $ R \subset R' $. 

Let us prove that $ L' = L $ and $ R' = R $. 
Let $ z = xgy \in L'GR' $. 
Since $ \cS(\nu) = \cS(\cK) $, we may take $ z_n \in \cS(\lambda_n) $ such that $ z_n \to z $. 
Since $ \cK $ is abelian, we have 
$ \lambda_n * \lambda_n^{-1} = \lambda_n^{-1} * \lambda_n = \eta $, 
and by Proposition \ref{prop: mu*nu LR} we have 
$ \lambda_n^{L'} = \eta^{L'} = \eta^L $ 
and $ \lambda_n^{R'} = \eta^{R'} = \eta^R $. 
Hence we obtain $ x_n := z_n^{L'} \in \cS(\lambda_n^{L'}) = \cS(\eta^L) = L $ 
and $ y_n := z_n^{R'} \in \cS(\lambda_n^{R'}) = \cS(\eta^R) = R $, 
and thus $ x = \lim x_n \in L $ and $ y = \lim y_n \in R $, 
which shows $ L' = L $ and $ R' = R $. 

(v) 
Let $ \lambda \in \cK $. For $ z = xgy \in \cS(\lambda ) \subset \cS(\nu) = LGR $, 
since $ RL \subset H $, we have 
\begin{align}
xgy \in LgHR \subset LHRxgyLHR \subset \cS(\eta) \cS(\lambda) \cS(\eta) 
\subset \cS(\eta * \lambda * \eta) = \cS(\lambda) . 
\label{}
\end{align}
Hence we have $ \cS(\lambda) = LG_{\lambda }R $ 
for $ G_{\lambda} := \bigcup \{ gH : z=xgy \in \cS(\lambda) \} \subset G $, 
and we also have $ G_{\lambda} = \bigcup \{ Hg : z=xgy \in \cS(\lambda) \} $ similarly. 
Note that $ G_{\lambda } H = H G_{\lambda } = G_{\lambda } $. 
Take $ g_{\lambda} \in G $ such that 
$ H g_{\lambda}^{-1} \subset G_{\lambda^{-1}} $. Then we obtain 
\begin{align}
L H g_{\lambda}^{-1} G_{\lambda} R 
\subset L G_{\lambda^{-1}} RL G_{\lambda} R 
\subset \cS(\lambda^{-1}) \cS(\lambda) 
\subset \cS(\lambda^{-1} * \lambda) 
\subset \cS(\eta) = LHR , 
\label{eq: LHglam-1}
\end{align}
which yields that $ H g_{\lambda}^{-1} G_{\lambda} \subset H $ 
and hence $ G_{\lambda} = g_{\lambda} H $. 
Similarly, we obtain $ G_{\lambda} = H g_{\lambda} $. 

For any $ h \in H $ and $ g \in G \subset \cS(\nu) = \cS(\cK) $, we may take 
$ z_n = x_ng_ny_n \in \cS(\lambda_n) $ such that 
$ z_n \to g $ and consequently $ g_n \to g $. 
In a similar way to \eqref{eq: LHglam-1}, we have 
\begin{align}
g_n h g_n^{-1} \in (g_n H) (H g_n^{-1}) 
= G_{\lambda_n} G_{\lambda_n^{-1}} 
\subset \cS(\eta) = LHR , 
\label{}
\end{align}
which shows $ g_n h g_n^{-1} \in eLHRe = H $. 
Letting $ n \to \infty $, we obtain $ ghg^{-1} \in H $, 
which shows that $ H $ is a normal subgroup of $ G $. 
Since $ G $ and $ H $ are both compact, 
we see by \cite[Theorem 5.22]{HR1} that 
the quotient group $ G/H = \{ gH : g \in G \} $ is also compact. 
Let $ \pi:G \to G/H $ denote the natural projection. 

Since 
\begin{align}
\cS(\eta^R * \lambda * \eta^L) 
= \bar{ \cS(\eta^R) \cS(\lambda) \cS(\eta^L) } 
= \bar{ R L G_{\lambda} R L } \subset \bar{ H g_{\lambda} H H }  = g_{\lambda} H , 
\label{}
\end{align}
we obtain the convolution factorization 
\begin{align}
\lambda = \eta \lambda \eta 
= \eta^L * \omega_H * (\eta^R * \lambda * \eta^L) * \omega_H * \eta^R 
= \eta^L * \omega_{g_{\lambda} H} * \eta^R . 
\label{eq: conv factr lambda}
\end{align}

We now define the mapping $ F:\cK \to G/H $ by $ F(\lambda) := g_{\lambda} H $. 
For $ \lambda_1,\lambda_2 \in \cK $, then 
\begin{align}
\lambda_1 * \lambda_2 
= \eta^L * \omega_{g_{\lambda_1} H} * (\eta^R * \eta^L) * \omega_{g_{\lambda_2} H} * \eta^R 
= \eta^L * \omega_{(g_{\lambda_1} g_{\lambda_2} H)} * \eta^R , 
\label{}
\end{align}
since $ RL \subset H $, which shows that $ F $ is a group homomorphism. 
Injectivity of $ F $ is obvious by \eqref{eq: conv factr lambda}. 
Let $ g \in G $. As we have seen it above, we may take 
$ z_n = x_ng_ny_n \in \cS(\lambda_n) $ 
such that $ g_n \to g $ and $ g_n H = g_{\lambda_n} H $. 
Then, by \eqref{eq: conv factr lambda}, we have 
\begin{align}
\lambda_n = \eta^L * \omega_{g_{\lambda_n} H} * \eta^R 
\to \eta^L * \omega_{gH} * \eta^R =: \lambda . 
\label{}
\end{align}
This shows that $ \lambda \in \cK $ and $ F(\lambda) = gH $, 
which yields surjectivity of $ F $. 
Suppose $ \cK \ni \lambda_n \to \lambda \in \cK $.  
By \eqref{eq: conv factr lambda}, we have 
\begin{align}
\omega_{F(\lambda_n)} = \delta_e * \lambda_n * \delta_e 
\to \delta_e * \lambda * \delta_e = \omega_{F(\lambda)} 
\quad \text{in $ \cP(G) $}, 
\label{}
\end{align}
which shows by the continuity of the natural projection $ \pi $ that 
\begin{align}
\delta_{F(\lambda_n)} = \omega_{F(\lambda_n)} \circ \pi^{-1} 
\to \omega_{F(\lambda)} \circ \pi^{-1} = \delta_{F(\lambda)} 
\quad \text{in $ \cP(G/H) $}, 
\label{}
\end{align}
which implies $ F(\lambda_n) \to F(\lambda) $ and we have seen continuity of $ F $. 
Since $ \cK $ is compact and $ G/H $ is Hausdorff, 
we see by \cite[Theorem 9 of Chapter 5]{Kel} that $ F $ is homeomorphic. 
Since $ F(\mu * \eta) \in G/H $, 
we may take $ \gamma \in G $ such that $ F(\mu * \eta) = \gamma H $, 
and then we obtain \eqref{eq: conv factr muketa} 
since $ (\mu * \eta)^k = \mu^k = \mu^k * \eta $ and $ F $ is a group homomorphism. 

Finally, let us prove the representations \eqref{eq: cK G repre}. 
Since any $ \lambda \in \cK $ 
can be represented as $ \lambda = \lambda * \eta = \lim \mu^{n(k)} * \eta $, 
we see that $ \cK = \bar{ \{ \eta,\mu*\eta,\mu^2*\eta,\ldots \} } $. 
Since for any $ g \in G $ we have $ F(\lambda) = gH $ 
for some $ \lambda = \lim \mu^{n(k)} * \eta \in \cK $, we obtain 
$ gH = F(\lambda) = \lim F(\mu^{n(k)} * \eta) = \lim \gamma^{n(k)} H $ in $ G/H $, 
which yields $ G/H = \bar{ \{ H,\gamma H,\gamma^2 H,\ldots \} } $. 
}


\begin{thebibliography}{10}

\bibitem{MR2050898}
G.~Budzban and A.~Mukherjea.
\newblock Subsemigroups of completely simple semigroups and weak convergence of
  convolution products of probability measures.
\newblock {\em Semigroup Forum}, 68(3):400--410, 2004.

\bibitem{MR119041}
G.~Choquet and J.~Deny.
\newblock Sur l'\'{e}quation de convolution {$\mu =\mu \ast \sigma $}.
\newblock {\em C. R. Acad. Sci. Paris}, 250:799--801, 1960.

\bibitem{MR269770}
S.~T.~L. Choy.
\newblock Idempotent measures on compact semigroups.
\newblock {\em Proc. London Math. Soc. (3)}, 20:717--733, 1970.

\bibitem{Coh}
D.~L. Cohn.
\newblock {\em Measure theory}.
\newblock Birkh\"{a}user Advanced Texts: Basler Lehrb\"{u}cher. [Birkh\"{a}user
  Advanced Texts: Basel Textbooks]. Birkh\"{a}user/Springer, New York, second
  edition, 2013.

\bibitem{MR137789}
H.~S. Collins.
\newblock Convergence of convolution iterates of measures.
\newblock {\em Duke Math. J.}, 29:259--264, 1962.

\bibitem{MR136679}
H.~S. Collins.
\newblock Idempotent measures on compact semigroups.
\newblock {\em Proc. Amer. Math. Soc.}, 13:442--446, 1962.

\bibitem{MR205306}
I.~Csisz\'{a}r.
\newblock On infinite products of random elements and infinite convolutions of
  probability distributions on locally compact groups.
\newblock {\em Z. Wahrscheinlichkeitstheorie und Verw. Gebiete}, 5:279--295,
  1966.

\bibitem{MR876261}
P.~L. Davies and D.~N. Shanbhag.
\newblock A generalization of a theorem of {D}eny with applications in
  characterization theory.
\newblock {\em Quart. J. Math. Oxford Ser. (2)}, 38(149):13--34, 1987.

\bibitem{MR826359}
Y.~Derriennic.
\newblock Sur le th\'{e}or\`eme de point fixe de {B}runel et le th\'{e}or\`eme
  de {C}hoquet-{D}eny.
\newblock {\em Ann. Sci. Univ. Clermont-Ferrand II Probab. Appl.},
  (4):107--111, 1985.

\bibitem{MR276400}
J.~Duncan.
\newblock Primitive idempotent measures on compact semigroups.
\newblock {\em Proc. Edinburgh Math. Soc. (2)}, 17:95--103, 1970.

\bibitem{MR83681}
R.~Ellis.
\newblock A note on the continuity of the inverse.
\newblock {\em Proc. Amer. Math. Soc.}, 8:372--373, 1957.

\bibitem{MR1039473}
H.~Furstenberg and Y.~Katznelson.
\newblock Idempotents in compact semigroups and {R}amsey theory.
\newblock {\em Israel J. Math.}, 68(3):257--270, 1989.

\bibitem{MR108690}
I.~Glicksberg.
\newblock Convolution semigroups of measures.
\newblock {\em Pacific J. Math.}, 9:51--67, 1959.

\bibitem{MR169971}
M.~Heble and M.~Rosenblatt.
\newblock Idempotent measures on a compact topological semigroup.
\newblock {\em Proc. Amer. Math. Soc.}, 14:177--184, 1963.

\bibitem{HR1}
E.~Hewitt and K.~A. Ross.
\newblock {\em Abstract harmonic analysis. {V}ol. {I}}, volume 115 of {\em
  Grundlehren der Mathematischen Wissenschaften [Fundamental Principles of
  Mathematical Sciences]}.
\newblock Springer-Verlag, Berlin-New York, second edition, 1979.
\newblock Structure of topological groups, integration theory, group
  representations.

\bibitem{HM}
G.~H{\"o}gn{\"a}s and A.~Mukherjea.
\newblock {\em Probability measures on semigroups}.
\newblock Probability and its Applications (New York). Springer, New York,
  second edition, 2011.
\newblock Convolution products, random walks, and random matrices.

\bibitem{ISY}
Y.~Ito, T.~Sera, and K.~Yano.
\newblock Resolution of sigma-fields for multiparticle finite-state action
  evolutions with infinite past.
\newblock Preprint, arXiv:2008.12407.

\bibitem{Kal}
O.~Kallenberg.
\newblock {\em Foundations of modern probability}.
\newblock Probability and its Applications (New York). Springer-Verlag, New
  York, second edition, 2002.

\bibitem{MR3462}
Y.~Kawada and K.~It\^{o}.
\newblock On the probability distribution on a compact group. {I}.
\newblock {\em Proc. Phys.-Math. Soc. Japan (3)}, 22:977--998, 1940.

\bibitem{Kechris}
A.~S. Kechris.
\newblock {\em Classical descriptive set theory}, volume 156 of {\em Graduate
  Texts in Mathematics}.
\newblock Springer-Verlag, New York, 1995.

\bibitem{Kel}
J.~L. Kelley.
\newblock {\em General topology}.
\newblock Graduate Texts in Mathematics, No. 27. Springer-Verlag, New
  York-Berlin, 1975.
\newblock Reprint of the 1955 edition [Van Nostrand, Toronto, Ont.].

\bibitem{MR123348}
B.~M. Kloss.
\newblock Probability distributions on bicompact topological groups.
\newblock {\em Theor. Probability Appl.}, 4:237--270, 1959.

\bibitem{MR1083341}
K.-S. Lau and W.~B. Zeng.
\newblock The convolution equation of {C}hoquet and {D}eny on semigroups.
\newblock {\em Studia Math.}, 97(2):115--135, 1990.

\bibitem{MR199306}
Y.-F. Lin.
\newblock Not necessarily abelian convolution semigroups of probability
  measures.
\newblock {\em Math. Z.}, 91:300--307, 1966.

\bibitem{MR1135263}
C.-C. Lo and A.~Mukherjea.
\newblock Convergence in distribution of products of {$d\times d$} random
  matrices.
\newblock {\em J. Math. Anal. Appl.}, 162(1):71--91, 1991.

\bibitem{MR293687}
A.~Mukherjea.
\newblock On the convolution equation {$P=PQ$} of {C}hoquet and {D}eny for
  probability measures on semigroups.
\newblock {\em Proc. Amer. Math. Soc.}, 32:457--463, 1972.

\bibitem{MR423458}
A.~Mukherjea.
\newblock Limit theorems for convolution iterates of a probability measure on
  completely simple or compact semigroups.
\newblock {\em Trans. Amer. Math. Soc.}, 225:355--370, 1977.

\bibitem{MR556681}
A.~Mukherjea.
\newblock Limit theorems: stochastic matrices, ergodic {M}arkov chains, and
  measures on semigroups.
\newblock In {\em Probabilistic analysis and related topics, {V}ol. 2}, pages
  143--203. Academic Press, New York-London, 1979.

\bibitem{Muk}
A.~Mukherjea.
\newblock {\em Topics in products of random matrices}, volume~87 of {\em Tata
  Institute of Fundamental Research Lectures on Mathematics}.
\newblock Published by Narosa Publishing House, New Delhi; for the Tata
  Institute of Fundamental Research, Mumbai, 2000.

\bibitem{MR516743}
A.~Mukherjea and T.~C. Sun.
\newblock Convergence of products of independent random variables with values
  in a discrete semigroup.
\newblock {\em Z. Wahrsch. Verw. Gebiete}, 46(2):227--236, 1978/79.

\bibitem{MR296207}
A.~Mukherjea and N.~A. Tserpes.
\newblock Idempotent measures on locally compact semigroups.
\newblock {\em Proc. Amer. Math. Soc.}, 29:143--150, 1971.

\bibitem{MT}
A.~Mukherjea and N.~A. Tserpes.
\newblock {\em Measures on topological semigroups: convolution products and
  random walks}.
\newblock Lecture Notes in Mathematics, Vol. 547. Springer-Verlag, Berlin-New
  York, 1976.

\bibitem{MR543577}
A.~Nakassis.
\newblock Limit behavior of the convolution iterates of a probability measure
  on a semigroup of matrices.
\newblock {\em J. Math. Anal. Appl.}, 70(2):337--347, 1979.

\bibitem{Par}
K.~R. Parthasarathy.
\newblock {\em Probability measures on metric spaces}.
\newblock AMS Chelsea Publishing, Providence, RI, 2005.
\newblock Reprint of the 1967 original.

\bibitem{MR801345}
H.~Pfister.
\newblock Continuity of the inverse.
\newblock {\em Proc. Amer. Math. Soc.}, 95(2):312--314, 1985.

\bibitem{MR148793}
J.~S. Pym.
\newblock Idempotent measures on semigroups.
\newblock {\em Pacific J. Math.}, 12:685--698, 1962.

\bibitem{MR1010824}
C.~R. Rao and D.~N. Shanbhag.
\newblock Further extensions of the {C}hoquet-{D}eny and {D}eny theorems with
  applications in characterization theory.
\newblock {\em Quart. J. Math. Oxford Ser. (2)}, 40(159):333--350, 1989.

\bibitem{MR890374}
A.~Raugi.
\newblock Un th\'{e}or\`eme de {C}hoquet-{D}eny pour les semi-groupes
  ab\'{e}liens.
\newblock In {\em Th\'{e}orie du potentiel ({O}rsay, 1983)}, volume 1096 of
  {\em Lecture Notes in Math.}, pages 502--520. Springer, Berlin, 1984.

\bibitem{MR0118773}
M.~Rosenblatt.
\newblock Limits of convolution sequences of measures on a compact topological
  semigroup.
\newblock {\em J. Math. Mech.}, 9:293--305, 1960.

\bibitem{MR185636}
M.~Rosenblatt.
\newblock Products of independent identically distributed stochastic matrices.
\newblock {\em J. Math. Anal. Appl.}, 11:1--10, 1965.

\bibitem{Ros}
M.~Rosenblatt.
\newblock {\em Markov processes. {S}tructure and asymptotic behavior}.
\newblock Die Grundlehren der mathematischen Wissenschaften, Band 184.
  Springer-Verlag, New York-Heidelberg, 1971.

\bibitem{MR841094}
M.~Rosenblatt.
\newblock Convolution sequences of measures on the semigroup of stochastic
  matrices.
\newblock In {\em Random matrices and their applications ({B}runswick, {M}aine,
  1984)}, volume~50 of {\em Contemp. Math.}, pages 215--220. Amer. Math. Soc.,
  Providence, RI, 1986.

\bibitem{MR169969}
\v{S}. Schwarz.
\newblock Convolution semigroup of measures on compact noncommutative
  semigroups.
\newblock {\em Czechoslovak Math. J.}, 14(89):95--115, 1964.

\bibitem{MR169970}
\v{S}. Schwarz.
\newblock Product decomposition of idempotent measures on compact semigroups.
\newblock {\em Czechoslovak Math. J.}, 14(89):121--124, 1964.

\bibitem{MR114874}
K.~Stromberg.
\newblock Probabilities on a compact group.
\newblock {\em Trans. Amer. Math. Soc.}, 94:295--309, 1960.

\bibitem{MR272935}
T.-C. Sun and N.~A. Tserpes.
\newblock Idempotent measures on locally compact semigroups.
\newblock {\em Z. Wahrscheinlichkeitstheorie und Verw. Gebiete}, 15:273--278,
  1970.

\bibitem{MR1046339}
G.~J. Sz\'{e}kely and W.~B. Zeng.
\newblock The {C}hoquet-{D}eny convolution equation {$\mu=\mu*\sigma$} for
  probability measures on abelian semigroups.
\newblock {\em J. Theoret. Probab.}, 3(2):361--365, 1990.

\bibitem{MR175168}
A.~Tortrat.
\newblock Lois tendues, convergence en probabilit\'{e} et \'{e}quation {$P\ast
  P^{\prime} =P$}.
\newblock {\em C. R. Acad. Sci. Paris}, 258:3813--3816, 1964.

\bibitem{MR92921}
K.~Urbanik.
\newblock On the limiting probability distribution on a compact topological
  group.
\newblock {\em Fund. Math.}, 44:253--261, 1957.

\bibitem{MR679395}
J.~Wo\'{s}.
\newblock The convolution equation of {C}hoquet and {D}eny for probability
  measures on discrete semigroups.
\newblock {\em Colloq. Math.}, 47(1):143--148, 1982.

\bibitem{MR125901}
W.~\.{Z}elazko.
\newblock A theorem on {$B_{0}$} division algebras.
\newblock {\em Bull. Acad. Polon. Sci. S\'{e}r. Sci. Math. Astronom. Phys.},
  8:373--375, 1960.

\end{thebibliography}
\end{document}